
\input amstex
\magnification 1200
\documentstyle{amsppt}
\topmatter
\title
On the Variety of Rational Curves of Degree $d$ in $\Bbb P^n$.
\endtitle
\author P. I. Katsylo.
\endauthor
\address
Independend University of Moscow,
Bolshoi Vlasievskii, 11, Moscow, 121002, Russia.
\endaddress
\email katsylo\@katsylo.mccme.ru
\endemail
\date April 26, 1998
\enddate
\abstract
We construct the canonical structure of an irreducible
projective variety on the set of connected curves of
degree $d$ in $\Bbb P^n$ with rational components
(some components can be multiple). The set of rational
curves is open subset in the Zariski topology
under this structure.
\endabstract
\endtopmatter
\document
\define\C{\Bbb C}
\define\PP{\Bbb P}
\define\SS{\Bbb S}
\define\SL{\operatorname{SL}}
\define\GL{\operatorname{GL}}
\define\R{\Cal R}
\define\qqq{\hfill $\square$}

\bf \S0. \rm
Consider the set $\R_d^0$ of rational curves of degree
$d$ in $\PP^n$. The set $\R_d^0$ is a subset of the set
$\R_d$ of the connected curves of degree $d$ in $\PP^n$
with rational irreducible components ( some components
can be multiple). In the present article we construct
the structure of an irreducible projective variety on
the set $\R_d$. The subset $\R_d^0$ is open in $\R_d$
in the Zariski topology. The base of the construction
is the classical one developed by Caley, van der Waerden,
and Chow. The construction is nothing but an attaching
of the "Chow coordinates" and the "Caley form" to a
cycle in $\PP^n$. See the textbooks [2], [3] and [4]
for a classical discription and for a references
on the works of Caley, Van der Warden, and Chow.
In \S $1$ we give a "modernized" construction
(we consider the "Caley tensors" instead of the
"Caley forms" and so on). We believe that this approach
is more natural and simple. In \S $2$ we define the
canonical structure of an irreducible projective variety
on the set of connected curves of degre $d$ in $\Bbb P^n$
with rational components such that the subset of rational curves
is open in the Zariski topology.

\bf \S1. \rm
We use standard facts of the representation theory (see [1]).
The group $\GL_{n+1}$ acts canonically on the space $\C^{n+1}$ and
therefore, the canonical action of the group $\GL_{n+1}$ on the
space $\SS_\lambda (\C^{n+1})$ and
$\SS_\lambda ((\C^{n+1})^\lor)$ are canonically defined
for an arbitrary $\lambda$ (by $\SS_\lambda$ we denote
the Schur functor corresponding to the partition $\lambda$,
by $V^\lor$ we denote the space dual to the linear space $V$).
Consider the closed subvariety
$$
Z_k \subset \wedge^k(\C^{n+1})^\lor
$$
of decomposable $k$-forms and the grade algebra
$$
\C[Z_k] = \underset{d \ge 0} \to {\otimes} \C[Z_k]_d
$$
of regular functions on $Z_k$. Recall that there is
the $\GL_{n+1}$-moduli isomorphism
$$
\C[Z_k]_d \simeq \SS_{\mu(k,d)} (\C^{n+1}) \tag 1.1
$$
(where $\mu(k,d) = (d, \dots , d)$ ($k$ times)) such that
$s \in \SS_{\mu(k,d)} (\C^{n+1})$ correspond to the homogenious
function
(we denote it by the same letter)
$$
s(z) = \langle s \otimes \undersetbrace
{n \ times} \to {z \otimes \dots \otimes z} \rangle,
\ \ \ \ z \in Z
$$
(we denote by
$$
\langle \cdot \rangle : V^{\otimes m} \otimes (V^\lor)^{\otimes m}
\rightarrow \C
$$
the complete contraction). We identify
$\C[Z_k]_d$ and $\SS_{\mu(k,d)} (\C^{n+1})$ according
to the isomorphism (1.1).

The projectivization
$$
PZ_k \subset P(\wedge^k(\C^{n+1})^\lor)
$$
of the subvariety $Z_k$ is the image of the $k$-grassmanian
$Gr(k,(\C^{n+1})^\lor)$ under the Plucker embedding.
The divisors group on $PZ_k$ is generated by an hyperplane section.
Therefore, an arbitrary cycle of codimension
$1$ in $PZ_k$ has the form
$$
\{ \overline{z} \in PZ_k \ | \ s(z) = 0 \}
$$
for some $0 \neq s \in \SS_{\mu(k,d)}(\C^{n+1})$.

Let
$$
X = X_1 + \dots + X_m \subset \PP^n = P\C^{n+1}
$$
be a $(k-1)$-dimensional cycle in $\PP^n$.
Consider the cycle
$$
Y_X = Y_1 + \dots + Y_m,
$$
in $PZ_k$, where
$$
Y_i = \{ \overline{z} \in PZ_k \ | \ Ker(z) \cap X_i \neq \emptyset \}
\subset PZ_k
$$
(by $Ker(z)$ we denote the kernel of the form $z$). We have
$$
codim_{PZ_k}(Y_X) = 1
$$
and therefore,
$$
Y_X = \{ \overline{z} \in PZ_k \ | \ Ca_X(z) = 0 \}
$$
for some $0 \neq Ca_X \in \SS_{\mu(k,d)}(\C^{n+1})$.
It is easy to see that $d$ is the degree of the cycle $X$.
\proclaim{Definition 1.1} The tensor $Ca_X$ is called
the Caley tensor of the cycle $X$.
\endproclaim
Set
$$
\gathered
Ch_{k,d,n+1} = \{ \overline{s} \in P\SS_{\mu(k,d)}(\C^{n+1}) \ |
\ s \ is \ the \ Caley \ tensor \\
\ of \ some \ (k-1)-dimensional \ cycle
\ of \ degree \ d \ in \ \PP^n \}.
\endgathered
$$
\proclaim{Definition 1.2} The variety $Ch_{k,d,n+1}$ is called
the Chow variety.
\endproclaim
\proclaim{Theorem 1.3} $Ch_{k,d,n+1}$ is a closed
in $P\SS_{\mu(k,d)}(\C^{n+1})$ subvariety.
\endproclaim
\demo{Proof} See [2] or [3].
\enddemo
Consider some simplest cases of the Caley tensors and
the Chow varieties. We have for $d=1$:
$$
\SS_{\mu(k,1)}(\C^{n+1}) = \wedge^k\C^{n+1},
\ \ \ \ \ Ch_{k,1,n+1} = Gr(k,\C^{n+1})
$$
is a $k$-grassmanian. The cycle corresponding to the point
$s \in Gr(k,\C^{n+1})$ is a $(k-1)$-dimensional linear
subspace corresponding to $s$.

We  have for $k=n$:
$$
\SS_{\mu(k,d)}(\C^{n+1}) \simeq S^d(\C^{n+1})^\lor
$$
(as $\SL_{n+1}$-moduli),
$$
Ch_{n,d,n+1} = S^d(\C^{n+1})^\lor
$$
is the space of forms of degree $d$. $(n-1)$-dimensional cycle
corresponding to the point $s \in S^d(\C^{n+1})^\lor$ is
the zeros set of $s$.

\bf \S2. \rm
We suppose in this section that $k = 2$, i. e.
we consider curves in $\PP^n$.

Set
$$
\gathered
\gathered
\R_d^0 = \{ \overline{s} \in Ch_{2,d,n+1} \ |
\ s \ is \ the \ Caley \ tensor \\
\ of \ some \ rational \ curve \ of \ degree \ d
\ in \ \PP^n \},
\endgathered \\
\R_d = \overline{\R_d^0}
\endgathered
$$
(we get the closure here in $P\SS_{(d,d)}(\C^{n+1})$ in the
Zarisski topology),
$$
\partial \R_d = \R_d \setminus \R_d^0.
$$
It is clear that the point $s \in \partial \R_d$ correspond to
a reducible curve. From \S $1$ it follows that $\partial \R_d$ is
closed in $\R_d$ and therefore, $\R_d^0$ is open in $\R_d$
in the Zariski topology.
In this section we construct an uniformization
of the subset $\R_d^0 \subset \R_d$.

Let $e_0, \dots , e_n$ be the canonical basis in $\C^{n+1}$.
Consider the space $\C^2$
and let $z = (z_0, z_1)$ be the coordinates in $\C^2$
corresponding to the canonical basis of $\C^2$. The element
$$
f(z) = f_0(z) e_0 + \dots + f_n(z) e_n \in
S^d\C^{2 \lor} \otimes \C^{n+1}
$$
defines canonically the homogeneous of degree $d$ polynomial
mapping (we denote it by the same letter)
$$
f : \C^2 \rightarrow \C^{n+1}
\ \ \ \ \ z \mapsto f(z).
$$

An element $f \in S^d\C^{2 \lor} \otimes \C^{n+1}$ such that
$f_0, \dots , f_n$ have no common root defines the morphism
$$
\overline{f} : \PP^1 \rightarrow \PP^n,
\ \ \ \ \ \overline{z} \mapsto \overline{f(z)}.
$$
Let $U$ be the open in the Zarisski topology subset
of $S^d\C^{2 \lor} \otimes \C^{n+1}$ of  $f$ such that
$f_0, \dots , f_n$ have no common root and
$\overline{f}(\PP^1)$ is a curve of degree $d$
(in other words: degree of the mapping $\overline{f}$ on
its image is equal to $1$).

The group $\GL_2$ acts canonically on the space
$S^d\C^{2 \lor} \otimes \C^{n+1}$ and the subset
$U$ is $\GL_2$-invariant. It is easily proved that
the stabilizer $(\GL_2)_f$ of an element $f \in U$ coincides
with the kernel $Z$ of the action
$\GL_2 : S^d\C^{2 \lor} \otimes \C^{n+1}$.
Therefore, $U$ consists of $\GL_2$-orbits which are isomorphic
to $\GL_2/Z$.
Consider the set of $\GL_2$-orbits $U/\GL_2$ of $U$.
It follows from above that there is the canonical
structure of an analytical variety on the $U/\GL_2$
such that the canonical mapping of the sets
$U \rightarrow U/\GL_2$ is a morphism of analytical manifolds.
Consider the mapping of the sets
$$
\Gamma : U/\GL_2 \rightarrow \R_d,
\ \ \ \ \ \GL_2 \cdot f \mapsto \overline{Ca_{\overline{f}(\PP^1)}}.
$$
Since $\overline{f'}(\PP^1) = \overline{f''}(\PP^1)$
for $f', f'' \in U$ iff $f'$ and $f''$
lie on the same $\GL_2$-orbit, then $\Gamma$ maps
$U/\GL_2$ onto $\R_d^0$ biuniquely.
\proclaim{Theorem 2.1} There is
a homogeneous of degree $2d$
$\GL_2 \times \GL_{n+1}$-covariant
$$
\gamma : S^d \C^{2\lor} \otimes \C^{n+1} \rightarrow
\SS_{(d,d)}(\C^{n+1})
$$
such that
$$
\Gamma (\GL_2 \cdot f) = \overline{\gamma(f)}
$$
for $f \in U$.
\endproclaim
\demo{Proof}
Let $W$ be a $2$-dimensional linear space,
$w_1,w_2$ be a basis in $W$. Consider the following polyhomogeneous
$\GL_2 \times \GL_{n+1} \times \GL(W)$-covariants:
$$
\eta : (\C^{n+1})^\lor \otimes W \rightarrow
\wedge^2(\C^{n+1})^\lor
$$
("factorization" by $\GL(W)$, (the image of $\eta$ is the set
of decomposable $2$-forms),
$$
\beta : S^d \C^{2\lor} \otimes \C^{n+1} \times
(\C^{n+1})^\lor \otimes W \rightarrow S^d \C^{2\lor} \otimes W
$$
(tensor multiplication plus contruction),
$$
\gathered
\rho : S^d \C^{2\lor} \otimes W \rightarrow \C, \\
h_1 \otimes w_1 + h_2 \otimes w_2 \mapsto Res (h_1, h_2),
\endgathered
$$
($Res(h_1, h_2)$ is the resultant of the forms $h_1, h_2$),
$$
\gamma' = \rho \circ \beta : S^d \C^{2\lor} \otimes \C^{n+1}
\times (\C^{n+1})^\lor \otimes W \ \rightarrow \ \C.
$$
We have: $\gamma'(f,r) \equiv \gamma''(f,\eta(r))$, where
$$
\gamma'' : S^d \C^{2\lor} \otimes \C^{n+1}
\times \wedge^2(\C^{n+1})^\lor \ \rightarrow \ \C
$$
is a polyhomogeneous of polydegree $(2d, d)$
$\GL_2 \times \GL_{n+1}$-covariant.
Now define $\gamma$ in the following way:
$$
\gamma(f) = (\gamma''(f, \cdot) : Z_2
\rightarrow \C)
$$
(recall (see \S $1$) that we identify $\GL_{n+1}$-module
$\SS_{(d,d)}(\C^{n+1})$ and the space of homogeneous of degree $d$
regular functions on $Z_2 \subset \wedge^2(\C^{n+1})^\lor$). Consider
$r = v_1^* \otimes w_1 + v_2^*\otimes w_2 \in (\C^{n+1})^\lor \otimes W$
and the decomposable $2$-form $\eta(r)$. It is easy to see that
$\gamma''(f,\eta(r)) = \gamma'(f,r) = 0$ iff
the curve $\overline{f}(\PP^1)$ intersects
$$
Ker (\eta(r)) = \{ \overline{v} \in \PP^n |
\langle v, v_1^* \rangle = \langle v, v_2^* \rangle = 0 \}.
$$ \qqq
\enddemo
Therefore, $\Gamma$ is a morphism of analytical spaces,
$U/\GL_2$ is quasiprojective nonsingular variety,
$\Gamma$ maps biuniquely $U/\GL_2$ onto $\R_d^0$.

Let us discribe curves in $\PP^n$ lying on $\partial \R_d$.
\proclaim{Theorem 2.2} $Ca_X$ (where $X$ is a curve of degree
$d$ in $\PP^n$) lies on $\partial \R_d$
iff $X$ is a connected reducible curve with rational irreducuble
components.
\endproclaim
\demo{Proof}
It is clear that if $Ca_X \in \partial \R_d$, then $X$ is connected
reducible curve with rational irreducible components.

Now suppose that $X$ is a connected reducible curve
of degree $d$ in $\PP^n$ with rational irreducible
components.
Let us prove that $Ca_X \in \partial \R_d$. Fix an arbitrary
metric $dist(\cdot,\cdot)$ on $\PP^n$.
It is enough to prove that for every $\epsilon$ there is a rational
curve $X_\epsilon$ such that for every point
$x \in X$ there exists a point $x_\epsilon \in X_\epsilon$
such that $dist(x,x_\epsilon) \le \epsilon$.
By using induction we reduce the deal to the case then
$X$ consists of two irreducible components $X_1$ and $X_2$.
One can assume that $X_1 = \overline{f}(\PP^1)$, $deg(f)=d_1$,
$X_2 = \overline{g}(\PP^1)$, $deg(f)=d_2$,
$f(1,0) = g(0,1) = A = e_0 + \dots + e_n$.
Consider
$F_\epsilon(z) = f_0(z) g_0(\epsilon z_0, z_1) e_0 +
\dots + f_n(z) g_n(\epsilon z_0, z_1) e_n$. For a fixed
$(z_0,z_1)$, $z_1 \neq 0$ we have: $\overline{F_\epsilon(z)}
\rightarrow \overline{f(z)} \in X_1$ as $\epsilon \rightarrow 0$.
For $z_0 = 1$, $z_1 \rightarrow 0$, $\epsilon \rightarrow 0$,
$z_1/\epsilon = t$ we have: $\overline{F_\epsilon(z)}
\rightarrow \overline{g_0(1,t)e_0 + \dots + g_n(1,t)e_n} \in X_2$
as $\epsilon \rightarrow 0$. From this it follows the statement.
\enddemo

\enddocument